\documentclass[11pt]{article}
\usepackage{dsfont}
\usepackage{multirow}
\usepackage{mathrsfs}
\usepackage{amsmath}
\usepackage{mathtools}
\usepackage{cases}
\usepackage{amssymb}
\usepackage{graphicx}
\usepackage{pstricks}
\usepackage{pstricks-add}
\usepackage{epic}
\usepackage{color}
\usepackage{calc}
\usepackage{tikz}
\usepackage{paralist}
\usepackage{cite}
\usepackage[noblocks]{authblk}

\newtheorem{definition}{Definition}[section]
\newtheorem{theorem}[definition]{Theorem}
\newtheorem{lemma}[definition]{Lemma}

\newtheorem{conjecture}[definition]{Conjecture}

\newtheorem{proposition}[definition]{Proposition}

\numberwithin{equation}{section}

\def\Ga{\Gamma}

\def\la{\lambda}
\def\ge{\geq}
\def\le{\leq}
\def\tt{\theta}
\def\witi{\widetilde}

\begin{document}
\title{A spectral characterization of the $s$-clique extension of the triangular graphs}
\author[a]{Ying-Ying Tan}
\author[b,c]{Jack H. Koolen\footnote{Corresponding author. Email addresses: koolen@ustc.edu.cn (J.H. Koolen), tansusan1@ahjzu.edu.cn (Y.-Y. Tan), xzj@mail.ustc.edu.cn (Z.-J. Xia). \\
 }}
\author[d]{Zheng-Jiang Xia}
\affil[a]{\small School of Mathematics $\&$ Physics, Anhui Jianzhu University, Hefei, Anhui, 230201, PR China.}
\affil[b]{School of Mathematical Sciences, University of Science and Technology of
China, Hefei, Anhui, 230026, PR China}
\affil[c]{Wen-Tsun Wu Key Laboratory of the CAS, School of Mathematical Sciences, University of Science and Technology of
China, Hefei, Anhui, 230026, PR China}
 \affil[d]{School of Finance, Anhui University of Finance and Economics,
 Bengbu, Anhui, 233030, PR China.}

\maketitle

\noindent {\bf Abstract:} A regular graph is co-edge regular if there exists a constant $\mu$ such that any two distinct and non-adjacent vertices have exactly $\mu$ common neighbors. In this paper, we show that for
integers $s\ge 2$,and $n$ large enough, any co-edge-regular graph which is cospectral with the $s$-clique extension of the triangular graph $T(n)$ is exactly the $s$-clique extension of the triangular graph $T(n)$.


\medskip

\large{This paper is dedicated to the memory of Prof. Slobadan Simi$\rm \acute{c}$.}

\section{Introduction}
All graphs in this paper are simple and undirected. For definitions related to distance-regular graphs, see\cite{bcn89,vkt16}.
Before we state the main result, we give more definitions.

Let $G$ be a simple connected graph on vertex set $V(G)$, edge set $E(G)$ and adjacency matrix $A$. The eigenvalues of $G$ are the eigenvalues of $A$. Let $\lambda_0,\lambda_1,\ldots,\lambda_t$ be the distinct eigenvalues of $G$ and $m_i$ be the multiplicity of $\lambda_i$ ($i=0,1,\ldots,t$). Then the multiset
$\{\lambda_0^{m_0},\lambda_1^{m_1},\ldots,\lambda_t^{m_t}\}$ is called the \emph{spectrum} of $G$. Two graphs are called \emph{cospectral} if they have the same spectrum. Note that a graph $H$ cospectral with a $k$-regular graph $G$ is also $k$-regular.

Recall that a regular graph is co-edge-regular if there exists a constant $\mu$ such that any two distinct and non-adjacent vertices have exactly $\mu$ common neighbors. Our main result in this paper is as follows:
\begin{theorem}\label{thm1}
Let $\Gamma$ be a co-edge-regular graph with spectrum
$$\{(2sn-3s-1)^1,(sn-3s-1)^{n-1},(-s-1)^{\frac{n^2-3n}{2}},(-1)^{\frac{(s-1)n(n-1)}{2}}\},$$
where $s\geq2$ and $n\geq1$ are integers. If $n\geq 48s$, then $\Gamma$ is the $s$-clique extension of the triangular graph $T(n)$.
\end{theorem}

This paper is a follow-up paper of Hayat, Koolen and Riaz\cite{hayat}. They showed a similar result
 for the square grid graphs. In that paper, they gave the
following conjecture.

\begin{conjecture}\cite{hayat}
Let $\Ga$ be a connected co-edge-regular graph with four distinct eigenvalues. Let $t\ge 2$ be an integer and $|V(\Ga)|=n(\Ga)$. Then there exists a constant $n_t$ such that, if $\theta_{min}(\Ga)\ge -t$
and $n(\Ga)\ge n_t$ both hold, then $\Ga$ is the $s$-clique extension of a strongly regular graph for some
$2\le s\le t-1$.
\end{conjecture}

This conjecture is wrong as the $p\times q$-grids $(p >q\ge 2)$ show. So we would like to modify this
conjecture as follows.
\begin{conjecture}
Let $\Ga$ be a connected co-edge-regular graph with parameters $(n, k, \mu)$ having four distinct eigenvalues. Let $t\ge 2$ be an integer. Then there exists a constant $n_t$ such that,
if $\theta_{min}(\Ga)\ge -t$, $n\ge n_t$ and $k< n-2-\frac{(t-1)^2}{4}$, then either $\Ga$ is the
$s$-clique extension of a strongly regular graph for  $2\le s\le t-1$ or $\Ga$
is a $p\times q$-grid with $p>q\ge 2$.
\end{conjecture}

The reason for the valency condition is, that in \cite{yang}, it was shown that for $\la\ge 2$, there
exist  constants $C(\la)$ such that a connected $k$-regular co-edge-regular graph with order $v$
and smallest eigenvalue at least $-\la$ satisfies one of the following:
\begin{itemize}
\item [(i)] $v-k-1\le \frac{(\la-1)^2}{4}+1$, or;
\item [(ii)]Every pair of distinct non-adjacent vertices have at most $C(\la)$ common neighbours.
\end{itemize}

Koolen et al.\cite{koolen}
improved this result by showing that one can take $C(\la)=(\la-1)\la^2$ if $k$ is much larger
than $\la$.
This paper is part of the project to show the conjecture for $t=3$.

Another motivation comes from the lecture notes \cite{terw}. In these notes, Terwilliger shows that
any local graph of a thin $Q$-polynomial distance-regular graph is co-edge-regular and has at most
five distinct eigenvalues. So it is interesting to study co-edge-regular graphs with a few distinct
eigenvalues.

We mainly follow the method of Hayat et al.\cite{hayat}. The main difference is that we simplify the method of Hayat
et al. when we  show that every vertex lies on exactly two lines.
This leads to a better bound for which we can show this. This will also improve the bound given in
the result of Hayat et al.

\section{Preliminaries}
\subsection{Definitions}

For two distinct vertices $x$ and $y$, we write $x\sim y$ (resp. $x \nsim y$) if they are adjacent (resp. nonadjacent) to each other. For a vertex $x$ of $G$, we define $N_G(x)=\{y\in V(G)\mid y\sim x \}$, and $N_G(x)$ is called the neighborhood of $x$. The graph induced by $N_G(x)$ is called the \emph{local graph} of $G$ with respect to $x$ and is denoted by $G(x)$. We denote the number of common neighbors between two distinct vertices $x$ and $y$ by $\lambda_{x,y}$ (resp. $\mu_{x,y}$) if $x\sim y$ (resp. $x\nsim y$).

A graph is called \emph{regular} if every vertex has the same valency.
A regular graph $G$ with $n$ vertices and valency $k$ is called \emph{co-edge-regular} with parameters $(n, k, \mu)$
if any two nonadjacent vertices have exactly $\mu=\mu(G)$ common neighbors. In addition, if any two adjacent vertices
have precisely $\la=\la(G)$ common neighbors,  then $G$ is called \emph{strongly regular} with
parameters $(n, k, \la, \mu)$.
A graph $G$ is called \emph{walk-regular} if the number of closed walks of length $r$ from a given vertex $x$ is
independent of the choice of $x$ for all $r$, that is to say, for any $x$, $A^{r}_{xx}$ is constant for all $r$, where $A$ is the adjacency matrix of $G$.

Let $X$ be a set of size $t$. The \emph{Johnson graph} $J(t, d)$ $(t\ge 2d)$ is a graph with vertex set ${X \choose d}$, the set
of $d$-subsets of $X$, where two $d$-subsets are adjacent whenever they have $d-1$ elements in common. $J(t, 2)$ is
the \emph{triangular graph} $T(t)$. Recall that a \emph{clique} (or a complete graph) is a graph in which every pair of vertices is adjacent. A \emph{coclique} is a graph that any two distinct vertices are nonadjacent.  A \emph{$t$-clique} is a clique with $t$-vertices and is denoted by $K_t$. The line graph of $K_t$ is  also the triangular graph $T(t)$ which is strongly regular with parameters $({t \choose 2},2t-4,t-2,4)$ and spectrum
$\{(2t-4)^1,(t-4)^{t-1},(-2)^{\frac{t^2-3t}{2}}\}$.

The Kronecker product $M_1\otimes M_2$ of two matrices $M_1$ and $M_2$ is obtained by replacing the $ij$-entry of $M_1$ by $(M_1)_{ij}M_2$ for all $i$ and $j$. Note that if $\tau$ and $\eta$ are eigenvalues of $M_1$ and $M_2$ respectively, then $\tau\eta$ is an eigenvalue  of $M_1\otimes M_2$.

\subsection{Interlacing}
\begin{lemma}[\cite{hae}, Interlacing]\label{int}
Let $N$ be a real symmetric $n\times n$ matrix with eigenvalues $\theta_1\ge \ldots \ge \theta_n$ and $R$ be a real $n\times m$ $(m < n)$
matrix with $R^TR=I$. Set $M=R^TNR$ with eigenvalues $\mu_1\ge\ldots \ge \mu_m$. Then,
\begin{enumerate}
\item  the eigenvalues of $M$ interlace those of $N$, i.e.
$$\tt_i\ge \mu_i\ge \tt_{n-m+i}, ~~~i=1, 2,\ldots, m,$$
\item if the interlacing is tight, that is, there exists an integer $j\in \{1,2,\ldots, m\}$ such that\\ $\tt_i=\mu_i$ for
$1\le i\le j$ and $\tt_{n-m+i}=\mu_i$ for $j+1\le i\le m$, then $RM=NR.$
\end{enumerate}
\end{lemma}

In the case that $R$ is permutation-similar to
$\left(
  \begin{array}{cc}
    I & O \\
    O & O \\
  \end{array}
\right)$, then $M$ is just a principal submatrix of $N$.

Let $\pi=\{V_1,\ldots, V_m\}$ be the partition of the index set of the columns of $N$ and let $N$ be partitioned according to $\pi$ as\\

$\left(
  \begin{array}{ccc}
    N_{1,1} &\ldots & N_{1,m} \\
    \vdots & \ddots & \vdots \\
    N_{m,1} & \ldots & N_{m,m} \\
  \end{array}
\right)$,\\

 where $N_{i, j}$ denotes the block matrix of $N$ formed by rows in $V_i$ and columns in $V_j$. The \emph{characteristic matrix} $P$ is the $n\times m$ matrix whose $j$th column is the characteristic vector of $V_j$ $(j=1, \ldots, m)$. \emph{The quotient matrix} of $N$ with respect to $\pi$ is the $m\times m$ matrix $Q$ whose entries are the average row sum of the blocks $N_{ij}$ of $N$, i.e.,
$$Q_{i,j}=\frac{1}{V_i}(P^TNP)_{i,j}.$$
The partition $\pi$ is called \emph{equitable} if each block $N_{i,j}$ of $N$ has constant row (and column) sum, i.e., $PQ=NP$. The following lemma can be shown by using Lemma \ref{int}.

\begin{lemma}\emph{\cite{god}}
Let $N$ be a real symmetric matrix with $\pi$ as a partition of the index set of its columns. Suppose $Q$ is the quotient matrix of $N$ with respect to $\pi$, then the following hold:
\begin{enumerate}
\item The eigenvalue of $Q$ interlace the eigenvalues of $N$.
\item If the interlacing is tight (as defined in Lemma \emph{\ref{int}(ii))}, then the partition $\pi$ is equitable.

\end{enumerate}
\end{lemma}

By an equitable partition of a graph, we always mean an equitable partition of its adjacency matrix $A$.

\subsection{Clique extensions of $T(n)$}
In this subsection, we define $s$-clique extensions of graphs and we will give some specific results for the $s$-clique extension of  triangular graphs.

Recall a $s$-clique is a
clique with $s$ vertices, where $s$ is a positive integer. The \emph{$s$-clique extension} of a graph $G$ with $|V(G)|$ vertices is the graph $\widetilde{G}$ obtained from $G$ by replacing
each vertex $x\in V(G)$ by a clique $\widetilde{X}$ with $s$ vertices, satisfying $\widetilde{x}\sim \widetilde{y}$ in $\widetilde{G}$ if and only if $x\sim y$ in $G$, where $\widetilde{x}\in \widetilde{X}, \widetilde{y}\in \widetilde{Y}$.
If $\witi{G}$ is a $s$-clique extension of $G$, then the adjacency matrix of $\witi{G}$ is $(A+I_{|V(G)|})\otimes J_s-I_{s|V(G)|}$, where $J_s$ is the all-ones matrix
of size $s$ and $I_{|V(G)|}$ is the identity matrix of size $|V(G)|$. In particular, if
$G$ has $t+1$ distinct eigenvalues and its spectrum is
\begin{align}\label{eq2.1}
\theta_0^{m_0}, \theta_1^{m_1},\ldots, \theta_t^{m_t},
\end{align}
then the spectrum of $\witi{G}$ is
\begin{align}\label{eq2.2}
\{(s(\theta_0+1)-1)^{m_0}, (s(\theta_1+1)-1)^{m_1},\ldots, (s(\theta_t+1)-1)^{m_t}, (-1)^{(s-1)(m_0+m_1+\cdots+m_t)}\}.
\end{align}
Note that if the adjacency matrix $A$ of a connected regular graph $G$ with $|V(G)|$ vertices and valency $k$ has four distinct eigenvalues $\{\theta_0=k, \theta_1, \theta_2, \theta_3\}$,
then $A$ satisfies the following equation (see \cite{hoff}):
\begin{align}\label{eq2.3}
A^3-( \sum\limits_{i=1}^3 \theta_i )A^2+ ( \sum_{1\leq i<j\leq 3}\theta_i\theta_j ) A-\theta_1\theta_2\theta_3I=\frac{\prod_{i=1}^3(k-\theta_i)}{|V(G)|}J.
\end{align}
This implies that $G$ is walk-regular, see \cite{vdam}.

Now we assume $\Ga$ is a cospectral graph with the $s$-clique extension of the triangular graph $T(n)$, where $s\ge 2$ and $n\ge 4$ are  integers. Then
by (\ref{eq2.1}) and (\ref{eq2.2}), the graph $\Ga$ has spectrum
\begin{align}\label{eq2.4}
\{\tt_0^{m_0}, \tt_1^{m_1}, \tt_2^{m_2},\tt_3^{m_3}\}=\{(s(2n-3)-1)^1, (s(n-3)-1)^{n-1}, (-s-1)^{\frac{n^2-3n}{2}}, (-1)^{(s-1)\frac{n(n-1)}{2}}\}.
\end{align}
Note that $\Ga$ is regular with valency $k$, where $k=(s-1)+2(n-2)s=s(2n-3)-1$.
Using (\ref{eq2.3}), we obtain
$$A^3+(3+4s-sn)A^2+((3-n)s^2+(8-2n)s+3)A+(1-(n-4)s-(n-3)s^2)I=4s^2(2n-3)J.$$
Therefore,
\begin{equation}\label{eq2.5}
A_{xy}^3=
\begin{cases}
2s^2n^2-2s^2n-6sn-3s^2+9s+2, &\enskip \text{if} ~x=y;\\
9s^2n+2sn-15s^2-8s-3-(3+4s-sn)\la_{xy}, &\enskip \text{if} ~x\sim y;\\
8s^2n-12s^2-(3+4s-sn)\mu_{xy}, &\enskip \text{if} ~x\nsim y.
\end{cases}
\end{equation}

The following result is known as the \emph{Hoffman bound}.

\begin{lemma}\emph( Cf. {\cite[Theorem 3.5.2]{bh12})}\label{lem2.1}
Let $X$ be a $k$-regular graph  with least eigenvalue $\tau$. Let $\alpha(X)$ be the size of maximum coclique in $X$.
Then $$\alpha(X)\le \frac{|X|(-\tau)}{k-\tau}.$$ If equality holds, then each vertex not in a coclique of size $\alpha(X)$ has exactly
$-\tau$ neighbours in it.
\end{lemma}

 Applying Lemma \ref{lem2.1} to the complement of $\Ga$, we obtain the following lemma.
\begin{lemma}\label{lem2.2}
For any clique $C$ of $\Ga$ with order $c$, we have
\begin{align*}
c\le s(n-1).
\end{align*}
If equality holds, then every vertex
$x\in V(\Ga)\setminus V(C)$ has exactly $2s$ neighbors in $C$.
\end{lemma}


\section{Lines in $\Ga$}
Recall that $\Ga$ is a graph that is cospectral with the $s$-clique extension of the triangular graph $T(n)$, where $s\ge 2$ and $n\ge 1$ are  integers. This
implies that $\Ga$ is walk-regular. Now we assume that $\Ga$ is also co-edge-regular, i.e., there exist  precisely $\mu=\mu(\Ga)$ common
neighbors between any two distinct nonadjacent vertices of $\Ga$. Note that for $\Ga$, we have $\mu=4s$ from the spectrum of the $s$-clique extension of $T(n)$.

Fix a vertex, denoted by $\infty$ and let $\Ga(\infty)$ be the local graph of $\Ga$ at vertex $\infty$.
Let $V(\Ga(\infty))=\{x_1, x_2,\ldots, x_k\}$, where $k=s(2n-3)-1$. Let $x_i$ have valency $d_i$ inside $\Ga(\infty)$ for $i=1,2,\ldots, k$.
Because $\Ga$ is walk-regular, the number of closed walks through a fixed vertex $\infty$ of length $3$ and $4$ only depend
 on the spectrum. This means that the number of edges  in $\Ga(\infty)$ is determined by the spectrum and as $\Ga$ is co-edge-regular, we also see that the number of walks of length $2$ in $\Ga(\infty)$ is determined by the spectrum of $\Ga$.
 This implies these numbers are the same as in a local graph of the $s$-clique extension of $T(n)$.

Let $\Delta$ be the $s$-clique extension of $T(n)$.
Fix a vertex $u$ of $\Delta$. Then
there are $s-1$ vertices with valency $(s-2)+2s(n-2)$ and $2s(n-2)$ vertices with valency $s(n-2)+2(s-1)$
in the local graph of $T(n)$ with respect
to a fixed vertex.
Using (\ref{eq2.5}),
this implies that the sum of valencies and the sum of
square of valencies of vertices in $\Ga(\infty)$ are constant, and are given by the following equations.
\begin{align}\label{eq3.6}
\sum\limits_{i=1}^k d_i=2\varepsilon=2s^2n^2-2s^2n-6sn-3s^2+9s+2,
\end{align}
\begin{align}\label{eq3.7}
\sum\limits_{i=1}^k (d_i)^2=2sn(s^2n^2-6sn-6s^2+10s+8)+9s^3+3s^2-24s-4,
\end{align}
where $\varepsilon$ is the number
of edges inside $\Ga(\infty)$.
By (\ref{eq3.6}) and (\ref{eq3.7}), we obtain
\begin{align}\label{eq3.8}
\sum\limits_{i=1}^k (d_i-(sn-2))^2=(s-1)s^2(n-3)^2.
\end{align}

It turns out that (\ref{eq3.8}) is of crucial importance in proving our main result.
Now we show the following lemma that will be used later.

\begin{lemma}\label{nocolique}
Fix a vertex $\infty$ of $\Ga$ and let $\Ga(\infty)$ be the local graph of $\Ga$ at $\infty$. Define $E:=\{y\sim\infty~|~d_y > \frac{3}{4}s(n-1)\}$ and let $e:=|E|$.
Let $F:=\{y\sim \infty\mid d_y\le \frac{3}{4}s(n-1)\}$ and                                                                                                                                               $f:=|F|$.
If $n\ge 55$, then the following hold.
\begin{asparaenum}[(1)]
\item $f\le 16(s-1)$.
\item The subgraph of $\Ga$ induced on $E$ is not complete.
\item The subgraph of $\Ga$ induced on $E$ does not contain a coclique of order three.
\end{asparaenum}
\end{lemma}
\textbf{Proof.}  Note that $f=k-e$. As
$\frac{3}{4}s(n-1)+1\le \frac{3}{4}(sn-2)$, by (\ref{eq3.8}), we obtain
\begin{equation}
\begin{aligned}
(s-1)s^2(n-3)^2=\sum_{y\sim \infty}(d_y-(sn-2))^2
&\ge \sum_{y\in F}(d_y-(sn-2))^2\\
&\ge \sum_{y\in F}(\frac{1}{4}(sn-2))^2\\
&=\frac{1}{16}f(sn-2)^2\\
&\ge \frac{1}{16}f(sn-s)^2.
\end{aligned}
\end{equation}
So $$f\le 16(s-1), $$ which
implies $f< \frac{1}{2}(sn-2)$ if $n\ge 55$ (and $s\ge 2$).
This means $$e=k-f>sn.$$
By Lemma \ref{lem2.2}, we obtain that $e$ is greater than the
order of a maximum size clique and hence the subgraph induced on $E$ is not complete.

Now we show that $E$ does not contain a coclique of order three.
Suppose $X\subset E$ is a coclique in $\Ga(\infty)$ with vertices $\{x_1, x_2, x_3\}$. 
Define $A_i~(i=1, 2, 3)$ such that
$$A_i:=\{y\sim \infty\mid y \sim x_i, y\nsim x_j~  \mbox{for all}~ x_j \in X, j\neq i\}\cup \{x_i\}.$$
Since $\Ga$ is co-edge-regular, the vertices $x_i$ and $x_j~(i\neq j)$ have at most $4s-1$ common neighbours.
By the inclusion-exclusion principle, we have
$$\frac{3\times(\frac{3}{4}s(n-1)+1)-k}{3}\le 4s-1.$$
This gives $n<54$.
This shows the lemma.
\hfill$\blacksquare$\\

A maximal clique of $\Ga$ is called a \emph{line} if it contains more than $\frac{3}{4}s(n-1)$ vertices.
We show the existence of lines of $\Ga$ in the following.

\begin{proposition}\label{pps3.2}
If $n\ge 48s\ge 96$, then for every vertex $\infty$, there are exactly two lines through $\infty$, say $C_1$ and $C_2$.
Denote $m:=|V(C_1)\cap V(C_2)\setminus \{\infty\}|$ and $\ell:=k+1-|V(C_1)\cup V(C_2)|$.
Then $m\le 4s-1$ and $\ell\le 16(s-1)$.
\end{proposition}
\textbf{Proof.}
Fix a vertex $\infty$ of $\Ga$, let $E=\{y\sim\infty~|~d_y> \frac{3}{4}s(n-1)\}$.
By Lemma \ref{nocolique}, a maximum coclique in $E$ has order two as $n\ge 48s\ge 55$.
Let $x_1$, $x_2$ be distinct nonadjacent vertices in $E$ and let $y\in E$. Then $y$ has at least one neighbour in $\{x_1, x_2\}$.

Let $A_i:=\{y\in E \mid y \sim x_i, y\nsim x_j~\mbox{for}~~ j=1,2, j\neq i$\} for $i=1,2$. Then
the subgraph induced on $A_i$ is complete for $i=1, 2$. Let $C_i$ be a maximal clique  containing the vertex set $\{\infty\}\cup A_i$ for $i=1, 2$.
Note that $C_1\neq C_2$ as $x_1\nsim x_2$. Let $M:=V(C_1)\cap V(C_2)\setminus \{\infty\}$ and $L:=V(\Ga(\infty))\setminus(V(C_1)\cup V(C_2))$. Let $m=|M|$ and $\ell=|L|$.
By the co-edge-regularity of $\Ga$, we have $m\le 4s-1$.
Let $F:=\{y\sim \infty\mid d_y\le \frac{3}{4}s(n-1)\}$ and $f:=|F|$.
We have, by Lemma \ref{nocolique}, that $f\le 16(s-1)$.

Suppose $x\in E\setminus (V(C_1)\cup V(C_2))$.
Then $x$ has at least $(\frac{3}{4}s(n-1)-(4s-2)-16(s-1))/2$ neighbours in at least
one of $C_1$ and $C_2$. If $n\ge 48s\ge 96$, then this number is at least $4s$, which is a contradiction. Hence
$E\subseteq V(C_1)\cup V(C_2)$.
So, $L\subseteq F$ and hence $\ell\le f\le 16(s-1)$ by Lemma \ref{nocolique}.
This shows that $|V(C_1)|+|V(C_2)|\ge k-\ell\ge k-16(s-1)$. Assume $|V(C_1)|\ge |V(C_2)|$, then we see that
$$|V(C_2)|\ge k-16(s-1)-s(n-1)>\frac{3}{4}s(n-1),$$ as $n\ge 48s\ge 96$. This gives that there are exactly two lines  through $\infty$.
\hfill$\blacksquare$\\

Now we prove the following property for lines through a vertex.

\begin{lemma}\label{lem3.3}
Fix a vertex $\infty$  of $\Gamma$ and let $C_1$ and $C_2$ be the two lines through $\infty$ with respective orders $c_1$ and $c_2$.
Let $L:=V(\Ga(\infty))\setminus(V(C_1)\cup V(C_2))$ and $M:=V(C_1)\cap V(C_2)\setminus \{\infty\}$, and $\ell=|L|$, $m=|M|\ge 0$.
If $n\geq 48s\ge 96$, then $\ell+m=s-1$ and $$s(n-3)+1\le c_i\le s(n-1)$$ for $i=1, 2$.
\end{lemma}
\textbf{Proof.}
Let $Q:=V(C_1)\Delta V(C_2)$, where $\Delta$ means ``symmetric difference''. Then, by Lemma \ref{lem2.2},
$|Q|\le |V(C_1)|+|V(C_2)|\le 2s(n-1)$.

Note that $Q$ is the complement of $L\cup M$ inside $V(\Ga(\infty))$.

For $y\in M$, we have
\begin{equation}\label{ym}
2sn-19s\le k-1-\ell\leq d_y\leq k-1=2sn-3s-2,
\end{equation}
by Proposition \ref{pps3.2}.

Now let
$y\in L$. Then $y$ has at least $4s-1$ neighbors in each of $C_1$ and $C_2$.
Hence, by Proposition \ref{pps3.2},  we obtain
\begin{equation}\label{yl}
d_y\le 2\times(4s-1)+\ell-1\le 2(4s-1)+16(s-1)-1\le 24s.
\end{equation}
By (\ref{eq3.8}), we obtain
\begin{equation}
\begin{aligned}
(s-1)s^2(n-3)^2&=\sum_{y\sim \infty}(d_y-(sn-2))^2\\
&\ge \sum_{y\in L}(d_y-(sn-2))^2+\sum_{y\in M}(d_y-(sn-2))^2\\
&\ge \ell((sn-s)-24s)^2+m((2sn-19s)-sn)^2\\
&=\ell s^2(n-25)^2+m s^2(n-19)^2\\
&\ge (\ell+m)s^2(n-25)^2.
\end{aligned}
\end{equation}
So $$\ell+m\le \frac{(s-1)(n-3)^2}{(n-25)^2}<s$$ if $n\ge 48s.$
Hence
\begin{equation}\label{eq3.9}
\ell+m\le s-1.
\end{equation}
This gives for $y\in L\cup M$, using (\ref{ym}), (\ref{yl}) and $l\le s-1$, that
$$d_y-(sn-2)\le k-1-(sn-2)=sn-3s.$$

Note that by (\ref{eq3.9}),
\begin{equation}\label{eq3.16}
\begin{aligned}
s(n-1)&\geq |V(C_j)|\geq1+k-s(n-1)-l\\
&\geq2sn-3s-s(n-1)-(s-1)\\
&=s(n-3)+1
\end{aligned}
\end{equation}
for $j=1,2$.

For $y\in V(\Ga(\infty))\setminus (L\cup M)$, we obtain
$$sn-4s\le |V(C_2)|-m-2\le d_y\le |V(C_2)|-1+4s-1+\ell\le sn+4s-3.$$
Hence $|d_y-(sn-2)|\le 4s.$

Now (\ref{eq3.8}) gives us
\begin{equation}
\begin{aligned}
(s-1)s^2(n-3)^2&=\sum_{y\sim \infty}(d_y-(sn-2))^2\\
&\le \sum_{y\in L\cup M}(d_y-(sn-2))^2+\sum_{y\in Q}(d_y-(sn-2))^2\\
&\le (\ell+m)s^2n^2+2s(n-1)(4s)^2.
\end{aligned}
\end{equation}
So $$\ell+m\ge \frac{(s-1)(n-3)^2-32s(n-1)}{n^2}>s-2,$$
if $n\ge 48s\ge 96.$ This implies $\ell+m=s-1$.
This shows the lemma.
\hfill$\blacksquare$\\

We obtain the following lemma immediately.

\begin{lemma}\label{lem3.5}
Fix a vertex $\infty$  of $\Gamma$ and let $C_1$ and $C_2$ be the two lines through $\infty$ with respective orders $c_1$ and $c_2$. Assume $m:=|V(C_1)\cap V(C_2)\setminus \{\infty\}|$. If $n\ge 48s$, then $c_1+c_2=2s(n-2)+2(m+1)$.
\end{lemma}
\textbf{Proof.}
Let $\ell=|V(\Ga(\infty))\setminus (V(C_1)\cup V(C_2))|$. Then we have
$$(c_1-m-1)+(c_2-m-1)+m+\ell=k=2sn-3s-1.$$
If $n\ge 48s$, then we have $\ell+m=s-1$ by Lemma~\ref{lem3.3}, hence
$$c_1+c_2=2s(n-2)+2(m+1).$$\hfill$\blacksquare$\\

In the next two sections, we will follow the method as used in Hayat et al.\cite{hayat}.

\section{The order of lines}
In this section, we will show the following lemma on the order of lines.

\begin{lemma}\label{lem4.1}
Let $s\geq2$ and $n\geq1$ be integers. Let $\Gamma$ be a co-edge-regular graph that is cospectral with the $s$-clique extension of the triangular graph $T(n)$. Let $q_i$ be the number of lines with order $s(n-3)+i$ for $i=1,\ldots,2s$ and $\delta=\sum_{i=1}^{2s}q_i$ be the number of lines in $\Gamma$. Assume $n\geq 48s$. Then
\begin{equation}\label{eq4.19}
\sum_{i=1}^{2s}(s(n-3)+i)q_i=sn(n-1)
\end{equation}
holds, and the number $\delta$ satisfies
\begin{equation}\label{eq4.20}
n\leq \delta\leq n+2.
\end{equation}
If $\delta=n$, then  $q_i=0$ for all $i<2s$, and $q_{2s}=n$.

\end{lemma}

\textbf{Proof.}
Assume $n\geq 48s$. By Proposition \ref{pps3.2}, any vertex of $\Ga$ lies on exactly two lines.
 Now consider the set
$$W=\{(x,C)|~x\in V(C),~\mbox{where}~ C ~\mbox{is a line}\}.$$
Then, by double counting, the cardinality of the set $W$, we see (\ref{eq4.19}).
Moreover, we see that
\begin{equation*}
\delta=\sum_{i=1}^{2s}q_i<\sum_{i=1}^{2s}\frac{s(n-3)+i}{s(n-3)}q_i=n+2+\frac{6}{n-3}.
\end{equation*}
Thus, if $n>10$, we obtain $$\delta\leq n+2.$$
On the other hand, we have
\begin{equation*}
\delta=\sum_{i=1}^{2s}q_i\geq\sum_{i=1}^{2s}\frac{s(n-3)+i}{s(n-1)}q_i=n.
\end{equation*}
This shows $\delta\geq n$, and $\delta=n$ implies that all lines have order $s(n-1)$, which means
$q_i\neq 0$ if and only if $i=2s$.
This completes the proof.
\hfill$\blacksquare$

\section{The neighborhood of a line}

In this section we will show the following proposition.

\begin{proposition}\label{pps5.1}
Let $\Gamma$ be a co-edge-regular graph that is cospectral with the $s$-clique extension of the triangular graph $T(n)$, where $s\geq2, n\geq1$ are integers.
If $n\geq 48s$, then $\Gamma$ contains exactly $n$ lines.
\end{proposition}

\textbf{Proof.}
In Lemma~\ref{lem4.1}, we have seen that the number $\delta$ of lines satisfies $n\leq \delta\leq n+2$. Now we assume that $n+1\leq\delta\leq n+2$, in order to obtain a contradiction.
Let $q_i$ be the number of lines of order $s(n-3)+i$ in $\Gamma$, where $i=1,\ldots,2s$. Let $h$ be minimal such that $q_h\neq 0$. Then clearly, $1\leq h\leq 2s$. Fix a line $C$ with exactly $s(n-3)+h$ vertices. Let $q'_i$ be the number of lines $C'$ with $s(n-3)+i$ vertices that intersect $C$ in at least one vertex. So $q_i\geq q'_i$. By Lemma~\ref{lem3.5}, we obtain
\begin{equation}\label{eq5.22}
|V(C)\cap V(C')|=\frac{h+i-2s}{2}
\end{equation}
By Proposition~\ref{pps3.2}, every vertex lies on exactly two lines, and hence we obtain
\begin{equation}\label{eq5.23}
\sum_{i=1}^{2s}q_i(\frac{h+i-2s}{2})\geq \sum_{i=1}^{2s}q'_i(\frac{h+i-2s}{2})=s(n-3)+h.
\end{equation}
Now multiply (\ref{eq5.23}) by 2 and subtract (\ref{eq4.19}) from obtained equation, we find
\begin{equation}
\begin{aligned}
\delta(h+s(1-n))&=\sum\limits_{i=1}^{2s} q_i(h+s(1-n))
&\ge s(-n^2+3n-6)+2h
\end{aligned}
\end{equation}
as $\delta=\sum\limits_{i=1}^{2s} q_i$.
This gives
\begin{equation*}
h(\delta-2)\geq 2s(n-3)+(\delta-n)s(n-1).
\end{equation*}
As $n+1\leq\delta\leq n+2$, we see
\begin{equation}\label{eq5.24}
hn\geq h(\delta-2)\geq 2s(n-3)+(\delta-n)s(n-1)\geq 2s(n-3)+s(n-1)=3sn-7s.
\end{equation}
Since $n\geq 48s$,  (\ref{eq5.24}) implies that $h\geq 3s$. This contradicts to $h\leq2s$.
This completes the proof.
\hfill$\blacksquare$

\section{Proof of the main result}
In this section we show our main result, Theorem~\ref{thm1}.

\textbf{Proof of Theorem~\ref{thm1}.} Assume $n\geq 48s$. By Propositions~\ref{pps3.2} and \ref{pps5.1} and Lemma~\ref{lem4.1}, we find that there are exactly $n$
lines, each of order $s(n-1)$, and every vertex $x$ in $\Gamma$ lies on exactly two lines. Moreover, by Lemma~\ref{lem3.5}, the two lines through any vertex $x$ have exactly $s$ vertices
in common, and every neighbor of $x$ lies in one of the two lines through $x$. Now consider the following equivalence relation $\mathcal{R}$ on the vertex set $V(\Gamma)$:
$x\mathcal{R}x'$ if and only if $\{x\}\cup N_{\Gamma}(x)=\{x'\}\cup N_{\Gamma}(x')$, where $x,x'\in V(\Gamma)$.

Every equivalence class under $\mathcal{R}$ contains $s$ vertices and it is the intersection of two lines. Let us define the graph $\hat{\Gamma}$ whose vertices are the equivalent classes and two classes, say $S_1$ and $S_2$, are adjacent in $\hat{\Gamma}$ if and only if any vertex in $S_1$ is adjacent to any vertex in $S_2$. Then $\hat{\Gamma}$ is a regular graph with valency $2n-4$, and $\Gamma$ is the $s$-clique extension of $\hat{\Gamma}$. Note that the spectrum of $\hat{\Gamma}$ is equal to
$$\{(2n-4)^1,(n-4)^{n-1},(-2)^{\frac{n^2-3n}{2}}\},$$
by the relation of the spectra of $\Gamma$ and $\hat{\Gamma}$, see~(\ref{eq2.1})~and~(\ref{eq2.2}). Since $\hat{\Gamma}$ is a connected regular graph with valency $2n-4$, and it has exactly three distinct eigenvalues, it follows that $\hat{\Gamma}$ is a strongly regular graph with parameters $(\binom{n}{2},2n-4,n-2,4)$.

As proved in \cite{cha59}, the triangular graphs are determined by the spectrum except when $n=8$, but we suppose $n$ is large enough then the assertion follows immediately.
This completes the proof.
\hfill$\blacksquare$

\section{Acknowledgements}
Jack Koolen is partially supported by the National Natural Science Foundation of China (No. 11471009 and No. 11671376) and Anhui Initiative in Quantum Information Technologies (No. AHY150000). Ying-Ying Tan is supported by the
National Natural Science Foundation of China (No. 11801007) and Natural Science Foundation of Anhui Province (No. 1808085MA17) and Doctoral Start-up foundation of Anhui Jianzhu University (No. 2018QD22). Zheng-jiang Xia is supported by the University Natural Science Research Project of Anhui Province (No.
KJ2018A0438).

\end{document}